# On $h$-basis

Li An-Ping

Beijing 100085, P.R. China
apli0001@sina.com

Abstract

A number set $\mathcal{A}$ is called a $h$-basis of the interval $[0, n]$, if each integer in $[0, n]$ may be represented a sum of at most $h$ elements of $\mathcal{A}$. To find the minimal $h$-basis is also called postage stamp problem. In this paper, we will present a new construction for $h$-basis, which is near the optimal one locally.

## 1. Introduction

Suppose that $\mathcal{A}$ and $\mathcal{B}$ are two number sets, as usual, define

$$\mathcal{A} + \mathcal{B} = \{a + b \mid a \in \mathcal{A}, b \in \mathcal{B}\}.$$

and for a non-negative integer $h$, $h\mathcal{A}$ represents the set of all possible sums of $h$ elements of $\mathcal{A}$, where the elements repetition is allowed, that is, $h\mathcal{A}$ is the summation of $h$'s sets $\mathcal{A}$. Suppose that $n$ is a natural number, a set $\mathcal{A}$ of non-negative integers is called an additive basis of order $h$ for the interval $[0, n]$, if $h\mathcal{A} \supseteq [0, n]$. Denoted by $n(h, \mathcal{A})$ the largest integer $n$ with that $h\mathcal{A} \supseteq [0, n]$, and for positive integer $k$, define

$$n(h, k) = \max_{|\mathcal{A}|=k} \{n(h, \mathcal{A})\}. \tag{1.1}$$

We call $\mathcal{A}$ is a extreme $h$-basis if $n(h, \mathcal{A}) = n(h, k), |\mathcal{A}| = k$. To find the precise value of $n(h, k)$ and design extreme $h$-basis is also called postage stamp problem.

Rohrbach [6] prove that

$$\left(\frac{k}{h}\right)^h \leq n(h, k) \leq \binom{k+h}{h}. \tag{1.2}$$

For the case $h = 2$, he presents a lower bound

$$n(2, k) \geq \frac{k^2}{4} + k + \delta, \qquad \delta \leq 1. \tag{1.3}$$

Hämmerer and Hofmeister [3] improve as

$$n(2, k) \geq \frac{10}{9} \frac{k^2}{4}. \tag{1.4}$$

And then, in papers [4], [5] it is further improved that

$$n(2, k) \geq \frac{2}{7} k^2. \tag{1.5}$$

It is not difficult to know that a construction of $h$-basis may be extended to a construction of $\lambda h$-basis, where $\lambda$ is an arbitrary positive integer. With this mechanism, and applying a unpublished result of R. Windecker for $n(3, k)$, Hofmeister [4] prove that

$$n(h, k) \geq (4/3)^{[h/3]} (8/7)^{[(h-3[h/3])/2]} (k/h)^h - O(h^{k-1}). \tag{1.6}$$

Correspondingly, for $h$ and $n$ given, we define $\zeta(h, n)$ as the minimal size in all the $h$-basis of interval $[0, n]$, namely

$$\zeta(h,n) = \min_{\mathcal{A}}\{|\mathcal{A}||n(h,\mathcal{A})\}. \tag{1.7}$$

Then (1.6) may be re-written as

$$\zeta(h,n) \leq \frac{n^{1/h}}{(4/3)^{1/3}}(h+o(h)). \tag{1.8}$$

For the detail of postage stamp problem, refer to see the survey paper of Alter and Barnett [1]. In this paper, we will present a new construction for $h$-basis, our main result is that

**Theorem 1.** If $n \geq \exp(h^2)$, then

$$\zeta(h,n) \leq n^{1/h}\left(h/e + 2.32 \log\log n + o(h)\right). \tag{1.9}$$

From Theorem 1, it directly follows that

**Corollary 1.** Let $\alpha_h = \sup_k\left\{\frac{n(h,k)^{1/h}}{k}\cdot\frac{h}{e}\right\}$, then $\lim_{h\to\infty}\alpha_h = 1$.

**2. The Proof of Theorem 1.**

In the new construction, we will apply two notions, $B_k$-sequence and the complement of a subset in finite group, which is introduced in paper [6] with a extension in this paper.

A set $\mathcal{A}$ is called $B_k$-sequence if all the sums of arbitrary $k$ elements of $\mathcal{A}$ are different. Denoted by $\Phi_k(n)$ the maximal size of $B_k$-sequences contained in the interval $[0,n]$. S.C. Bose and S. Chowla [2] provided a lower-bound of $\Phi_k(n)$

$$\Phi_k(n) \geq n^{1/k} + o(n^{1/k}). \tag{2.1}$$

Suppose that $G$ is an Abelian group of order $q$, $A$ is a non-empty subset of $G$, and $X_1, X_2, \ldots, X_k$ are $k$ subsets of $G$, we call the family $\{X_i\}_1^k$ is a $k$-complement of $A$ to $G$, if $G \subseteq A \oplus X_1 \oplus X_2 \oplus \cdots \oplus X_k$, $\left|\bigcup_{1\leq i\leq k} X_i\right|$ is the size of the complement, and denoted by $\partial^k(A,G)$ the minimal $k$-complement of $A$ to $G$. For convenience, write $\hat{C}_r^k = r(r+1)\cdots(r+k-1)/k!$, $m = |\partial^k(A,G)|$, then by the definition, it has $\hat{C}_m^k \cdot |A| \geq |G|$, and

$$\left|\partial^k(A,G)\right| \geq (2\pi k)^{1/2k} \frac{k}{e} \cdot \left(\frac{|G|}{|A|}\right)^{1/k} - \frac{k}{2}. \tag{2.2}$$

In the proof of Theorem 1 it will be employed a sharper upper-bound for $\left|\partial^k(A,G)\right|$.

**Lemma 1.** Suppose that $G$ is an Abelian group of order $q$, $A$ and $B$ are arbitrary two subsets of $G$, then for arbitrary integer $t \geq 0$, there is a subset $X$ of $t$ elements in $G$ such that

$$\left|B \setminus (A \oplus X)\right| \leq \left(1 - \frac{|A|}{q}\right)^t \cdot |B|. \tag{2.3}$$

*Proof.* We will take the induction on integer $t$, and at first prove (2.3) for $t = 1$.
For any $v \in G$, write $A_v = A \oplus v$, denoted by $\lambda = |A|$, then clearly

$$\bigcup_{v \in G} A_v = \lambda \otimes G,$$

where expression $\lambda \otimes G$ represents $\lambda$ copies of set $G$. So,

$$B \cap \bigcup_{v \in G} A_v = \bigcup_{v \in G} (B \cap A_v) = \lambda \otimes (B \cap G) = \lambda \otimes B.$$

Hence, $\sum_{v \in G} |B \cap A_v| \geq \lambda |B|$, and there is at least one $z \in G$, such that $|A_z \cap B| \geq \lambda |B|/|G|$, and,

$$\left|B \setminus (A \oplus v)\right| = |B| - |B \cap A_v| \leq \left(1 - \frac{\lambda}{q}\right) \cdot |B|.$$

This prove (2.3) for $t = 1$. Suppose that $X \subseteq G, |X| \leq t$, such that (2.3) hold. Let $B' = B \setminus (A \oplus X)$, replacing $B$ by $B'$, as the observation above for $t = 1$, there is $v \in G$, such that $|B' \setminus A_v| \leq (1 - \lambda/q) \cdot |B'|$, let $\tilde{X} = X \cup v$, then

$$\left|B \setminus (A \oplus \tilde{X})\right| = \left|B \setminus ((A \oplus X) \cup (A \oplus v))\right| = \left|(B \setminus (A \oplus X)) \setminus (A \oplus v)\right| = \left|B' \setminus (A \oplus v)\right|$$

$$\leq \left(1 - \frac{\lambda}{q}\right) \cdot |B'| \leq \left(1 - \frac{\lambda}{q}\right)^{t+1} \cdot |B|.$$

So the induction is finished. □

Lemma 1 is also applied in paper [6].

**Lemma 2.** Suppose that $G$ is an Abelian group of order $q \, (>1)$, $A$ is a non-empty subset of

$G$, let $k_0$ be the maximal integer $\kappa$ with that $\left((q\log q)/|A|\right)^{1/\kappa} \geq \log q + 2$, then

$$\left|\partial^k(A,G)\right| \leq \begin{cases} k \cdot \left\lceil \left((q\log q)/|A|\right)^{1/k} \right\rceil + \lceil \log q \rceil, & \text{if } k \leq k_0, \\ k_0 \cdot \left\lceil \left((q\log q)/|A|\right)^{1/k_0} \right\rceil + \lceil \log q \rceil, & \text{if } k > k_0. \end{cases} \quad (2.4)$$

*Proof.* At first, assume that $k \leq k_0$, and let $\alpha = |A|$, $t = \left\lceil \left((q\log q)/\alpha\right)^{1/k} \right\rceil$, by lemma 1, there is a subset $X_1$ with no more than $t$ elements, such that

$$|G \setminus (A+X_1)| \leq \left(1-\frac{\alpha}{q}\right)^t q \leq e^{-\alpha t/q} \cdot q \leq \left(1-\frac{\alpha t}{q}+\frac{\alpha^2 t^2}{2q^2}\right) q$$

i.e.

$$|(A+X_1)| \geq \left(\frac{\alpha t}{q}-\frac{(\alpha t)^2}{2q^2}\right) \cdot q = \left(1-\frac{(\alpha t)}{2q}\right)\alpha t$$

Let $A_1 = A + X_1$, and apply Lemma 1 to $A_1$ again, there is the subset $X_2$ with no more than $t$ elements, such that

$$|G \setminus (A_1+X_2)| \leq \left(1-\frac{|A_1|}{q}\right)^t q \leq e^{-|A_1|t/q} \cdot q \leq \left(1-\frac{|A_1|t}{q}+\frac{|A_1|^2 t^2}{2q^2}\right) q$$

i.e.

$$|A_1+X_2| \geq \left(\frac{|A_1|t}{q}-\frac{|A_1|^2 t^2}{2q^2}\right) q \geq \left(1-\frac{|A_1|t}{2q}\right)|A_1|t$$

$$\geq \left(1-\frac{|A_1|t}{2q}\right)\left(1-\frac{(\alpha t)}{2q}\right)\alpha t^2 \geq \left(1-\frac{\alpha t}{2q}-\frac{|A_1|t}{2q}\right)\alpha t^2$$

$$\geq \left(1-\frac{\alpha t}{2q}-\frac{\alpha t^2}{2q}\right)\alpha t^2. \quad (\because |A_1| \leq \alpha t)$$

By the induction, there are subsets $X_1, X_2, \cdots, X_{k-1}, |X_i| \leq t, 1 \leq i \leq k-1$, such that

$$|A+X_1+X_2+\cdots+X_{k-1}| \geq \left(1-\frac{\alpha t}{2q}-\frac{\alpha t^2}{2q}-\cdots-\frac{\alpha t^{k-1}}{2q}\right)\alpha t^{k-1}$$

$$\geq \left(1-\frac{\alpha t^k}{2q(t-1)}\right)\alpha t^{k-1}$$

Let $A_{k-1} = A + X_1 + X_2 + \cdots + X_{k-1}$, $t' = t + \lceil \log q \rceil$, apply Lemma 1 to $A_{k-1}$ again, there is a subset $X_k$ with no more than $t'$ elements, such that

$$|G \setminus (A_{k-1} + X_k)| \leq \left(1 - \frac{|A_{k-1}|}{q}\right)^{t'} \cdot q$$

$$\leq \left(1 - \frac{\alpha t^{k-1}}{q}\left(1 - \frac{\alpha t^k}{2q(t-1)}\right)\right)^{t'} \cdot q$$

$$\leq \left(1 - \frac{\log q}{t}\left(1 - \frac{\log q}{2(t-1)}\right)\right)^{t'} \cdot q$$

$$\leq \left(1 - \frac{\log q}{t'}\right)^{t'} \cdot q < 1.$$

So the first one of (2.4) is proved, and the second one is obvious for $|\partial^{k_1}(A, G)| \leq |\partial^{k_2}(A, G)|$, if $k_1 \geq k_2$. □

It seems that the second estimation of (2.4) for the case $k > k_0$ may be done better through further investigation, for which will not been used in this paper, so is simplified.

*The Proof of Theorem 1*:

Let $p = \lceil n^{1/h} \rceil$, $a, k$ be two integers, $1 \leq k \leq a < h$, which will be determined later. Denoted by $m = p^{h-a} \cdot (h-a)!$, then by (2.1), there is a $B_{h-a}$-sequence $\mathcal{B}$ in $[0, m]$ with $|\mathcal{B}| \doteq m^{1/(h-a)}$. Denoted by $\mathcal{H} = (h-a)\mathcal{B}$, it is clear that $|\mathcal{H}| = \hat{C}_{|\mathcal{B}|}^{(h-a)} \geq p^{h-a}$, where expression $\hat{C}_r^d = r(r+1)\cdots(r+d-1)/d!$. Let $q = p^{h-a+k}$, and $\mathcal{C} = \partial^k(\mathcal{H}, Z_q)$. By Lemma 2, it has

$$|\mathcal{C}| \leq k \times \left\lceil \left(\frac{q}{|\mathcal{H}|}\log q\right)^{1/k} \right\rceil + \lceil \log q \rceil \leq k \cdot p \cdot \lceil (\log n)^{1/k} \rceil + \lceil \log n \rceil.$$

Suppose that $\mathcal{A}$ is an extreme *h-basis* of $[0, hq]$, by (1.2),

$$|\mathcal{A}| \leq (p^{h-a+k} \cdot h)^{1/h} \cdot h$$
$$\leq p \times h^{1/h} \times h / n^{(a-k)/h^2}$$
$$\leq p \cdot \sigma. \qquad (\sigma = h^{1+(1/h)} / n^{(a-k)/h^2})$$

Let $\mathcal{D} = \{j \cdot p^i \mid 0 \leq j < p, h-a+k \leq i < h\}$, $\mathcal{G} = \mathcal{A} \cup \mathcal{B} \cup \mathcal{C} \cup \mathcal{D}$. We claim that set $\mathcal{G}$ is a $h$-basis of $[0, n]$: For any integer $z \in [0, n]$, $z \geq hq$, write $z = sq + r$, $0 \leq r < q$, $s \geq h$. Then, there are $x \in \mathcal{H}$, $y \in k\mathcal{C}$, such that $x + y \equiv r \mod q$, suppose that $x + y = tq + r$,

clearly, $t < h$, hence it has

$$z = x + y + (s-t)q \in (h-a)\mathcal{B} + k\mathcal{C} + (a-k)\mathcal{D}.$$

Let $\tau \ (\doteq 0.8415)$ be the zero of equation $e^t(1-t) - e^{-1} = 0$, and $k = \lceil (\log\log n)/\tau \rceil$, $a = k + \lceil 2\log h \rceil$, then it has

$$|\mathcal{G}| \leq |\mathcal{A}| + |\mathcal{B}| + |\mathcal{C}| + |\mathcal{D}|$$
$$\leq p \cdot \sigma + p((h-a)/e + o(h)) + k \cdot \lceil p \cdot (\log n)^{1/k} \rceil + \lceil \log n \rceil + p(a-k)$$
$$\leq p \times \left((h-a)/e + k \cdot (\log n)^{1/k} + (a-k) + o(h)\right)$$
$$\leq p \times \left((h-k)/e + k \cdot (\log n)^{1/k} + (1-(1/e))(a-k) + o(h)\right)$$
$$\leq p \times \left(h/e + 2.32 \times \log\log n + o(h)\right).$$

□

Corollary 1 indicates that the construction described above will near the optimal one locally when $h$ is larger and $\log\log n = o(h)$, however, due to the term $\log\log n$ it is impeded to more large interval of $n$, so, any one substantively improvement on Lemma 2 will bring corresponding improvement on Theorem 1.